\documentstyle[12pt]{article}

\title{{\bf ON ALTERNATIVE SUPERMATRIX REDUCTION}
}
\author{{\bf Steven Duplij} \thanks{Alexander von Humboldt Fellow}
\thanks{On leave of absence from
{\sl Theory Division, Nuclear Physics Laboratory,
Kharkov State University, KHARKOV 310077, Ukraine}}
\thanks{E-mail: duplij@physik.uni-kl.de}\\
{\sl Physics Department, University of Kaiserslautern,}\\
{\sl Postfach 3049, D-67653 KAISERSLAUTERN,}\\
{\sl Germany}}
\date{31 May, 1995
}
\input amssym.def
\input amssym
\newtheorem{definition}{Definition}
\newenvironment{remark}[1]{\vskip 5pt\noindent{\it Remark.\ }#1}{\vskip
5pt}
\newenvironment{example}[1]{\vskip 5pt\noindent{\it Example.\ }#1}{\vskip 5pt}
\newenvironment{proof}[1]{\noindent{\it Proof.\
}#1}{\hskip 3pt $\Box$\vskip 5pt}
\newtheorem{theorem}[definition]{Theorem}

\newtheorem{corollary}[definition]{Corollary}
\newtheorem{assertion}[definition]{Assertion}
\newtheorem{proposition}[definition]{Proposition}
\newcommand{\limfunc}{\mbox}
\begin{document}
\maketitle
     \begin{abstract}
We consider  a nonstandard odd reduction of
supermatrices (as compared with the standard even one) which arises
in connection with possible extension  of
manifold structure group reductions. The study was initiated by consideration
of the
generalized noninvertible superconformal-like transformations.
 The features of even- and odd-reduced
 supermatrices are investigated
on a par. They can be  unified into some kind of "sandwich" semigroups.
Also we define a special module over even- and odd-reduced
 supermatrix sets,
and the generalized Cayley-Hamilton theorem is proved for them.
 It is shown that the
odd-reduced supermatrices represent semigroup bands and Rees
matrix semigroups over a unit group.
\end{abstract}
\begin{flushright}KL-TH-95/15
\end{flushright}
\newpage
\section{Introduction}

According to the general theory of $G$-structures \cite{che1,gui,kobayashi}
 various
geometries are obtained by a reduction of a structure group of a manifold to
some subgroup $G$ of the tangent space endomorphisms. In the local approach
using coordinate description this means that one should reduce a
corresponding matrix in a given representation to some reduced form as a
matter of fact. In the most cases this form is triangle, because of the
simple observation from the ordinary matrix theory that the triangle
matrices preserve the shape and form a subgroup. In supersymmetric theories,
despite of appearance of odd subspaces and anticommuting variables, the
choice of the reduction shape remained the same\cite{howe,lot,ros/sch1,schw4},
and a ground reason of this was the fully identity of the supermatrix
multiplication with the ordinary one, and consequently the shape of the
matrices from a subgroup was the same. However in fine search of nontrivial
supersymmetric manifestations one can observe that the closure of
multiplication can be also achieved for other shapes, but due to existence
of zero divisors in the Grassmann algebra or in the ring over which a theory
is defined. So the meaning of the reduction itself could be extended
principally. Evidently, that some ''good'' properties of the transformations
could be lost in this direction, but opening of new possibilities, beauty
and interesting and unusual features which are distinctive for
supersymmetric case only are the sufficient price for the surprises arisen
and reason for them to investigate.

This paper was initiated by the study of superconformal symmetry semigroup
extensions \cite{dup6,dup7}. Indeed superconformal transformations \cite
{bar/fro/sch1,cra/rab,coh} appear as a result of the reduction of the
structure group matrix to the triangle form \cite{gid/nel3,gid/nel1}. Also,
the transition functions on semirigid surfaces \cite{dis/nel,gov/nel/won}
(see \cite{dol/ros/sch}) occurred in the description of topological
supergravity \cite{gov/nel/rey} have the same shape. In \cite{dup6} we
considered an alternative version of the reduction. The superconformal-like
transformations obtained in this way have many unusual features, e. g. they
are noninvertible and twist parity of the tangent space in the
supersymmetric basis\footnote{
This situation is different from the case of $Q$-manifolds \cite
{ale/kon/sch/zab}, where changing parity of the tangent space is done by
hand from the first definitions.
}.

Here we study the alternative reduction of supermatrices from a more
abstract viewpoint without connecting a special physical model.
\section{Preliminaries}

Let $\Lambda $ be\footnote{
The standard notations can be found in \cite{berezin,lei1}, and here we list
some of them needed only.
}
a commutative Banach ${\Bbb Z}_2$-graded superalgebra over a field ${\Bbb K}$
(where ${\Bbb K}={\Bbb R,}$ ${\Bbb C}$ or ${\Bbb Q}_p$) with a decomposition
into the direct sum: $\Lambda =\Lambda _0\oplus \Lambda _1$. The elements $a$
from $\Lambda _0$ and $\Lambda _1$ are homogeneous and have the fixed even
and odd parity defined as $\left| a\right| \stackrel{def}{=}\left\{ i\in
\left\{ 0,1\right\} ={\Bbb Z}_2|\,a\in \Lambda _i\right\} $. The even
homomorphism ${\frak m}_b:\Lambda \rightarrow {\Bbb B}$ , where ${\Bbb B}$ is a
purely even algebra over ${\Bbb K}$, is called a body map, if for any other
purely even algebra ${\Bbb A}$ and any homomorphism ${\frak m}_a:\Lambda
\rightarrow {\Bbb A}$ there is an even homomorphism ${\frak m}_{ab}:{\Bbb B}%
\rightarrow {\Bbb A}$ such that ${\frak m}_a={\frak m}_{ab}\circ {\frak m}_b$.
The
kernel of ${\frak m}_b$ is ${\Bbb S\ }\equiv \ker {\frak m}_b\stackrel{def}{=}%
\left\{ a\in \Lambda |\,{\frak m}\left( a\right) =0\right\} $ and is called
the soul sector of $\Lambda $. If there are exists an embedding ${\frak n}:%
{\Bbb B\hookrightarrow \Lambda }$ such that ${\frak m}\circ {\frak n}=id$, then
$%
\Lambda $ admits the body and soul decomposition $\Lambda ={\Bbb B}\oplus
{\Bbb S}$, and a soul map can be defined as ${\frak m}_s:\Lambda \rightarrow
{\Bbb S}$. Usually the isomorphism ${\Bbb B\ }\cong {\Bbb K}$ is implied (which
is not necessary in general and can lead to very nontrivial behavior of the
body). This is the case when $\Lambda $ is modeled with the Grassmann
algebras $\wedge \left( N\right) $ having $N$ generators \cite
{rog1,rab/cra1,vla/vol} or $\wedge \left( \infty \right) $ \cite
{rog4,boy/git,khr2}, or with the free graded-commutative Banach algebras $%
\wedge _BE$ over Banach spaces \cite{jad/pil,pes7,bry1}. The soul ${\Bbb S}$
is obviously a proper two-sided ideal of $\Lambda $ which is generated by $%
\Lambda _1$. In case $\Lambda $ is a Banach algebra (with a norm $\left|
\left| \cdot \right| \right| $) soul elements are quasinilpotent \cite{iva2}%
, which means $\forall a\in {\Bbb S},\lim \limits_{n\rightarrow \infty
}\left| \left| a\right| \right| ^{1/n}=0$. But in the infinite-dimensional
case quasinilpotency of the soul elements does not necessarily lead to their
nilpotency ($\forall a\in {\Bbb S}\,\exists n,\,a^n=0$) \cite{pes5}. These
facts allow us to consider noninvertible morphisms on a par with invertible
ones (in some sense), which gives, in proper conditions, many interesting
and nontrivial results (see \cite{dup6,dup10,dup11}).

The $\left( p|q\right) $-dimensional linear model superspace $\Lambda ^{p|q}$
over $\Lambda $ (in the sense of \cite{rog1,vla/vol,khr2,dewitt}) is the
even sector of the direct product $\Lambda ^{p|q}=\Lambda _0^p\times \Lambda
_1^q$. The even morphisms $\limfunc{{\rm Hom}}_0\left( \Lambda
^{p|q},\Lambda ^{m|n}\right) $ between superlinear spaces $\Lambda
^{p|q}\rightarrow \Lambda ^{m|n}$ are described by means of $\left(
m+n\right) \times \left( p+q\right) $-supermatrices (for details see \cite
{berezin,lei1}). In various physical applications supermatrices are reduced
to some suitable form which is necessary for concrete consideration. For
instance, in the theory of super Riemann surfaces \cite{fri,ros/sch/vor1}
the $\left( 1+1\right) \times \left( 1+1\right) $-supermatrices describing
holomorphic morphisms of the tangent bundle have a triangle shape \cite
{gid/nel1,gid/nel3}.

Here we consider a special alternative reduction of supermatrices and study
its features. We note that the supermatrix theory per se has many own
problems \cite{bac/fel1,ebn,hus/nie} and unexpected conclusions (e.g. the
lowering of the degree of characteristic polynomials comparing to the
standard Cayley-Hamilton theorem\cite{urr/mor,urr/mor1}).

For transparency and clarity we confine ourselves to $\left( 1+1\right)
\times \left( 1+1\right) $-supermatrices\footnote{which
will allow us not to melt ideas
by large formulas, and only this size will be used in the following
consideration} , and generalization to the $\left( m+n\right) \times \left(
p+q\right) $ case is straightforward and can be mostly done by means of
simple changing of notations.

\section{Structure of $\limfunc{{\rm Mat}} _\Lambda \left( 1|1\right) $}

In the standard basis in $\Lambda ^{1|1}$ the elements from $\limfunc{{\rm
Hom}}%
_0\left( \Lambda ^{1|1},\Lambda ^{1|1}\right) $ are described by the $\left(
1+1\right) \times \left( 1+1\right) $-supermatrices \cite{berezin}
\begin{equation}
\label{1}M\equiv \left(
\begin{array}{cc}
a & \alpha \\
\beta & b
\end{array}
\right) \in \limfunc{{\rm Mat}} _\Lambda \left( 1|1\right)
\end{equation}

\noindent where $a,b\in \Lambda _0,\,\alpha ,\beta \in \Lambda _1$ (in the
following we use Latin letters for elements from $\Lambda _0$ and Greek
letters for ones from $\Lambda _1$). For sets of matrices we also use
corresponding bold symbols, e. g. $\ {\bf M}\stackrel{def}{=}\left\{ M\in
\limfunc{{\rm Mat}} _\Lambda \left( 1|1\right) \right\} $. In this simple $%
\left( 1|1\right) $ case the supertrace defined as $\limfunc{str}:\limfunc{%
Mat}_\Lambda \left( 1|1\right) \rightarrow \Lambda _0$ and
Berezinian defined as $\limfunc{{\rm Ber}} :\limfunc{{\rm Mat}} _\Lambda \left(
1|1\right) \setminus \left\{ M|\,{\frak m}_b\left( b\right) =0\right\}
\rightarrow \Lambda _0$ are
\begin{equation}
\label{2}\limfunc{str}M=a-b,
\end{equation}
\begin{equation}
\label{3}\limfunc{{\rm Ber}} M=\frac ab+\frac{\beta \alpha }{b^2}.
\end{equation}

Now we define two types of possible reductions of $M$ on a par and study
some of their properties simultaneously.
\begin{definition}
{\sl Even-reduced supermatrices} are elements from
$\limfunc{{\rm Mat}}_\Lambda \left( 1|1\right) $ having the form
\begin{equation}
\label{4}S\equiv \left(
\begin{array}{cc}
a & \alpha \\
0 & b
\end{array}
\right) \in \limfunc{{\rm RMat}} _\Lambda ^{\,even}\left( 1|1\right) .
\end{equation}

\noindent {\sl Odd-reduced supermatrices} are elements from $\limfunc{{\rm
Mat}} %
_\Lambda \left( 1|1\right) $ having the form
\begin{equation}
\label{5}T\equiv \left(
\begin{array}{cc}
0 & \alpha \\
\beta & b
\end{array}
\right) \in \limfunc{{\rm RMat}} _\Lambda ^{\,odd}\left( 1|1\right) .
\end{equation}
\end{definition}
The name of the odd-reduced supermatrices follows naturally from
$\limfunc{{\rm Ber}}  T=\beta \alpha /b^2\Rightarrow \left(
\limfunc{{\rm Ber}}  T\right) ^2=0$ and
\begin{equation}
\label{5a}\limfunc{{\rm Ber}}  T^2=\limfunc{{\rm Ber}} \left(
\begin{array}{cc}
\alpha \beta & \alpha b \\
\beta b & b^2-\alpha \beta
\end{array}
\right) =0.
\end{equation}

The explanation of the ground of the notations $S$ and $T$ comes from the
fact that the even-reduced supermatrices give superconformal transformations
which describe morphisms of the tangent bundle over the super Riemann
surfaces \cite{gid/nel1}, while the odd-reduced supermatrices give the
superconformal-like transformations twisting the parity of the $\left(
1|1\right) $ tangent superspace in the standard basis (see
\cite{dup6,dup10}).

\begin{assertion}
${\bf M}$ is a direct sum of diagonal ${\bf D}$ and
anti-diagonal (secondary diagonal) ${\bf A}$ supermatrices (the even and odd
ones in the notations of \cite{berezin})
\begin{equation}
\label{8}{\bf M=D\oplus A},
\end{equation}

\noindent where
\begin{equation}
\label{9}D{\bf \equiv }\left(
\begin{array}{cc}
a & 0 \\
0 & b
\end{array}
\right) \in {\bf D}\equiv \limfunc{{\rm Mat}} _\Lambda ^{\,Diag}\left(
1|1\right) ,
\end{equation}
\begin{equation}
\label{10}A{\bf \equiv }\left(
\begin{array}{cc}
0 & \alpha \\
\beta & 0
\end{array}
\right) \in {\bf A}\equiv \limfunc{{\rm Mat}} _\Lambda ^{\,Adiag}\left(
1|1\right) ,
\end{equation}
\noindent and ${\bf D\subset S}$ and ${\bf A\subset T}$.
\end{assertion}
For the reduced supermatrices one finds
\begin{equation}
\label{st}{\bf S\cap T}=\left(
\begin{array}{cc}
0 & \alpha  \\
0 & b
\end{array}
\right) {\bf \neq \emptyset .}
\end{equation}

Nevertheless, the following observation explains the fundamental role of $%
{\bf S}$ and ${\bf T}$.

\begin{proposition}
The Berezians of even- and odd-reduced supermatrices are additive components
of the full Berezinian
\begin{equation}
\label{st1}\limfunc{{\rm Ber}}M=\limfunc{{\rm Ber}}S+\limfunc{{\rm Ber}}T.
\end{equation}
\end{proposition}

The first term in (\ref{st1}) covers all subgroups of even-reduced
supermatrices from $\limfunc{Mat}_\Lambda \left( 1|1\right) $, and only it
was considered in the applications. But the second term is dual to the first
in some sense and corresponds to all subsemigroups of odd-reduced
supermatrices from $\limfunc{Mat}_\Lambda \left( 1|1\right) $%
\footnote{
The relation (\ref{st1}) is a supersymmetric version of the obvious equality
$\det M=\det D+\det A$, when $D$ and $A$ from (\ref{8}) and
(\ref{9}) are
ordinary matrices. The problem is that
 for $A$ being a supermatrix $\limfunc{{\rm Ber}}\,A$ is
not defined at all.
}.

\section{Invertibility and ideals of $\limfunc{{\rm Mat}} _\Lambda \left(
1|1\right) $}

Denote the set of invertible elements of ${\bf M}$ by ${\bf M^{*}}$, and $%
{\bf I=M\setminus M^{*}}$. In \cite{berezin} it was proved that ${\bf M^{*}}%
=\left\{ M\in {\bf M}|\,{\frak m}_b\left( a\right) \neq 0\wedge {\frak m}%
_b\left( b\right) \neq 0\right\} $. Then similarly ${\bf S^{*}}=\left\{ S\in
{\bf S}|\,{\frak m}_b\left( a\right) \neq 0\wedge {\frak m}_b\left( b\right)
\neq 0\right\} $ and ${\bf T^{*}}=\emptyset $, i. e. the odd-reduced
matrices are noninvertible and ${\bf T}\subset {\bf I}$. Consider the
invertibility structure of $\limfunc{{\rm Mat}} _\Lambda \left( 1|1\right)
$ in more detail. Let us denote
\begin{equation}
\label{10a}
\begin{array}{ccc}
{\bf M^{\prime }} & = & \left\{ M\in
{\bf M|\,}{\frak m}_b\left( a\right) \neq 0\right\} , \\ {\bf M^{\prime \prime
}} & = & \left\{ M\in
{\bf M|\,}{\frak m}_b\left( b\right) \neq 0\right\} , \\ {\bf I^{\prime }} & =
& \left\{ M\in
{\bf M|\,}{\frak m}_b\left( a\right) =0\right\} , \\ {\bf I^{\prime \prime }}
& = & \left\{ M\in {\bf M|\,}{\frak m}_b\left( b\right) =0\right\} .
\end{array}
\end{equation}

\noindent Then ${\bf M}={\bf M^{\prime }\cup I^{\prime }=M^{\prime \prime
}\cup I^{\prime \prime }}$ and ${\bf M^{\prime }\cap I^{\prime }=\emptyset }$
, ${\bf M^{\prime \prime }\cap I^{\prime \prime }=\emptyset }$, therefore $%
{\bf M}^{*}={\bf M^{\prime }}\cap {\bf M^{\prime \prime }}$ and ${\bf T}%
\subset {\bf M}^{\prime \prime }$. The Berezinian $\limfunc{{\rm Ber}} M$ is
well-defined for the matrices from ${\bf M}^{\prime \prime }$ only and is
invertible when $M\in {\bf M}^{*}$, but for the matrices from ${\bf M}%
^{\prime }$ the inverse $\left( \limfunc{{\rm Ber}} M\right) ^{-1}$ is
well-defined
and is invertible when $M\in {\bf M}^{*}$ too \cite{berezin}.

Under the ordinary matrix multiplication the set ${\bf M}$ is a semigroup of
all $\left( 1|1\right) $ supermatrices \cite{mca1,put2,put1},
and the set ${\bf M}^{*}$ is a
subgroup of ${\bf M}$. In the standard basis ${\bf M}^{*}$ represents the
general linear group $GL_\Lambda \left( 1|1\right) $ \cite{berezin}.
According to the general definitions \cite{cli/pre1} a subset ${\bf I}%
\subset {\bf M}$ is called a right (left) ideal of the semigroup ${\bf M}$,
if ${\bf I}\cdot {\bf M}\subset {\bf I}$ (${\bf M}\cdot {\bf I}\subset {\bf I%
}$), where the point denotes the standard matrix set multiplication: ${\bf %
A\cdot B}\stackrel{def}{=}\left\{ \bigcup AB|\,A\in {\bf A},\,B\in {\bf B}%
\right\} $. An isolated ideal satisfies the relation \cite{cli/pre1}
\begin{equation}
\label{11}M_1M_2\in {\bf I}\Rightarrow M_1\in {\bf I}\vee M_2\in {\bf I,}
\end{equation}

\noindent and  a filter ${\bf F}$ of the semigroup ${\bf M}$ is
defined by
\begin{equation}
\label{12}M_1M_2\in {\bf F}\Rightarrow M_1\in {\bf F\wedge }M_2\in {\bf F}.
\end{equation}

\begin{proposition}

1) The sets ${\bf I}$, ${\bf I^{\prime }}$ and ${\bf I^{\prime }}$ are
isolated ideals of ${\bf M}$.

2) The sets ${\bf M}^{*}$, ${\bf M}^{^{\prime }}$ and ${\bf M}^{^{\prime
\prime }}$ are filters of the semigroup ${\bf M}$.

3) The sets ${\bf M^{\prime }}$ and ${\bf M^{\prime \prime }}$ are not
subgroups\footnote{as it was incorrectly translated in \cite
{berezin}, pp. 95, 103 (in the original Russian edition, Moscow, 1983, pp.
89, 93, the sets ${\bf M^{\prime }}$ and ${\bf M^{\prime \prime }}$, denoted
as $G^{\prime }Mat\left( 1,1|\Lambda \right) $ and $G^{\prime \prime
}Mat\left( 1,1|\Lambda \right) $ respectively, are called semigroups).
},
but subsemigroups of ${\bf M}$, which are ${\bf M^{\prime }=M}^{*}\cup {\bf J%
}^{\prime }$ and ${\bf M^{\prime \prime }=M}^{*}\cup {\bf J^{\prime \prime }}
$ with the isolated ideals ${\bf J}^{\prime }={\bf M^{\prime }}\setminus
{\bf M}^{*}={\bf M^{\prime }}\cap {\bf I^{\prime \prime }}$ and ${\bf J}%
^{\prime \prime }={\bf M^{\prime \prime }}\setminus {\bf M}^{*}={\bf %
M^{\prime \prime }}\cap {\bf I^{\prime }}$ respectively.

4) The ideal of the semigroup ${\bf M}$ is
\begin{equation}
\label{13}{\bf I}={\bf I}^{\prime }\cup {\bf J}^{\prime }={\bf I}^{\prime
\prime }\cup {\bf J}^{\prime \prime }.
\end{equation}
\end{proposition}

\begin{proof}
Let $M_3=M_1M_2$, then $a_3=a_1a_2+\alpha _1\beta _2$ and $%
b_3=b_1b_2+\beta _1\alpha _2$. Taking the body part we derive
\begin{equation}
\label{14}{\frak m}_b\left( a_3\right) ={\frak m}_b\left( a_1\right) {\frak m}%
_b\left( a_2\right) ,
\end{equation}
\begin{equation}
\label{15}{\frak m}_b\left( b_3\right) ={\frak m}_b\left( b_1\right) {\frak m}%
_b\left( b_2\right) .
\end{equation}

1) The left-hand side of (\ref{14}) and (\ref{15}) vanishes iff the first or
second multiplier of the right-hand side equals zero. Then use (\ref{11}).

2) The left-hand side of (\ref{14}) and (\ref{15}) does not vanish iff the
first and second multiplier of the right-hand side does not equal zero. Then
use (\ref{12}).

3) ${\bf J}^{\prime }$ consists of the matrices with ${\frak m}_b\left(
a\right) \neq 0$, but ${\frak m}_b\left( b\right) =0$, and ${\bf J^{\prime
\prime }}$ consists of the matrices with ${\frak m}_b\left( a\right) =0$, but
$%
{\frak m}_b\left( b\right) \neq 0$.

4) The ideal of ${\bf M}$ consists of the matrices with ${\frak m}_b\left(
a\right) =0$ or ${\frak m}_b\left( b\right) =0$. Then use 3) and the
definitions.
\end{proof}
\begin{remark}
Since the ideals ${\bf J^{\prime }}$ and ${\bf %
J^{\prime \prime }}$ are isolated, i. e. ${\bf J^{\prime }\cap M}^{*}={\bf %
J^{\prime }\cap M}^{*}=\emptyset $, they cannot be represented as sequences
of elements from the group ${\bf M}^{*}$ (viz. no one noninvertible element
can be derived from a sequence of invertible ones, see, e. g. \cite{iva2}),
and so the statement ''that any element of $G^{\prime \prime }Mat\left(
p,q|\Lambda \right) $ (the semigroup ${\bf M}^{\prime \prime }$ here, and so
the notation $G^{\prime \prime }...$ misleads) is the limit of a sequence
from $GMat\left( p,q|\Lambda \right) $ (the group ${\bf M}^{*}$ here)''
holds only for invertible elements from ${\bf M}^{\prime \prime }$, i.e.
belonging ${\bf M^{*}}$, and it means that elements from ${\bf M}^{*}$ can
be obtained from a sequence of elements from ${\bf M}^{*}$, which is simply
a group action.
\end{remark}
\begin{assertion}
For the odd-reduced matrices from (\ref{10a}) it follows ${\bf T}\subset
{\bf I}^{\prime }$ and ${\bf A}\cap {\bf M^{\prime }=A}\cap {\bf M^{\prime
\prime }=\emptyset }$.
\end{assertion}

\section{Multiplication properties of odd-reduced supermatrices}

In general the odd-reduced matrices do not form a semigroup, since
\begin{equation}
\label{15a}T_1T_2=\left(
\begin{array}{cc}
\alpha _1\beta _2 & \alpha _1b_2 \\
b_1\beta _2 & b_1b_2+\beta _1\alpha _2
\end{array}
\right) \neq T.
\end{equation}

But from (\ref{15a}) it follows that
\begin{equation}
\label{15b}
\begin{array}{ccc}
{\bf T}\cdot {\bf T}\cap {\bf T}\neq \emptyset & \Rightarrow & \alpha \beta
=0, \\
{\bf T}\cdot {\bf T}\cap {\bf S}\neq \emptyset & \Rightarrow & \beta b=0,
\end{array}
\end{equation}

\noindent which can take place, because of the existence of zero divisors in
$\Lambda $.

\begin{proposition}

1) The subset ${\bf T}^{SG}\subset {\bf T}$ of the odd-reduced matrices
satisfying $\alpha \beta =0$ form an odd-reduced subsemigroup of ${\bf M}$.

2) In the odd-reduced semigroup ${\bf T}^{SG}$ the subset of matrices with $%
\beta =0$ is a left ideal, and one with $\alpha =0$ is a right ideal, the
matrices with $b=0$ form a two-sided ideal.

\end{proposition}

\subsection{Semigroup band representations}

Let
\begin{equation}
\label{b0}Z_\alpha \left( t\right) =\left(
\begin{array}{cc}
0 & \alpha t \\
\alpha  & 1
\end{array}
\right) \in {\bf Z}_\alpha {\bf \subset T}^{SG},
\end{equation}
i.e. ${\bf Z}_\alpha $ is a set of the odd-reduced matrices parameterized by
the even parameter $t \in \Lambda_0$. Then ${\bf Z}_\alpha $ is a semigroup
under the
matrix multiplication ($\alpha $ numbers the semigroups) which is isomorphic
to a one parameter semigroup with the multiplication
\begin{equation}
\label{b1}\left\{ t_1\right\} *_\alpha \left\{ t_2\right\} =\left\{
t_1\right\} .
\end{equation}
This semigroup is called a right zero semigroup ${\cal Z}_R=\left\{ \bigcup
\left\{ t\right\} ;*_\alpha \right\} $ and plays an important role (together
with the left zero semigroup ${\cal Z}_L$ defined in a dual manner) in the
general semigroup theory (e.g., see \cite{cli/pre1}, Theorem 1.27, and \cite
{howie}).

Let
\begin{equation}
\label{b11}B_\alpha \left( t,u\right) =\left(
\begin{array}{cc}
0 & \alpha t \\
\alpha u & 1
\end{array}
\right) \in {\bf B}_\alpha {\bf \subset T}^{SG},
\end{equation}
then ${\bf B}_\alpha $ is a matrix semigroup (numbered by $\alpha $) which
is isomorphic to a two $\Lambda_0$-parametric semigroup ${\cal B}=\left\{
\bigcup
\left\{ t,u\right\} ;*_\alpha \right\} $, where the multiplication is
\begin{equation}
\label{b2}\left\{ t_1,u_1\right\} *_\alpha \left\{ t_2,u_2\right\} =\left\{
t_1,u_2\right\} .
\end{equation}
Here every element is an idempotent (as in the previous case too), and so
this is a rectangular band multiplication \cite{howie,petrich2}.

Let $C_\alpha \left( t,u,v\right) =\left(
\begin{array}{cc}
0 & \alpha t \\
\alpha u & v
\end{array}
\right) \in {\bf C}_\alpha {\bf \subset T}^{SG}$, then ${\bf C}_\alpha $ is
a matrix semigroup isomorphic to a semigroup ${\cal B}_G=\left\{ \bigcup
\left\{ t,u,v\right\} ;*_\alpha \right\} $ where the multiplication is
\begin{equation}
\label{b3}\left\{ t_1,u_1,v_1\right\} *_\alpha \left\{ t_2,u_2,v_2\right\}
=\left\{ t_1v_2,u_2v_1,v_1v_2\right\} .
\end{equation}
The parameter $v$ describes the difference of an element from an idempotent,
since $\left\{ t,u,v\right\} ^2-\left\{ t,u,v\right\} =\left\{ t\left(
v-1\right) ,u\left( v-1\right) ,v\left( v-1\right) \right\} $.

\begin{assertion}
The one and two parametric subsemigroups of the semigroup of odd-reduced
supermatrices ${\bf T}^{SG}$ having vanishing Berezinian represent semigroup
bands, viz. the left and right zero semigroups and rectangular bands.
\end{assertion}
\begin{theorem}
The continuous supermatrix representation of the Rees matrix semigroup over
a unit group $G=e$ (see \cite{cli/pre1,howie}) is given by formulas (\ref{b0}%
) and (\ref{b11}).
\end{theorem}

\subsection{''Square root'' of even-reduced supermatrices}

Consider the second equation in (\ref{15b}).
\begin{proposition}
The elements $T^{\sqrt{S}}$from the subset ${\bf T}^{\sqrt{S}}\subset {\bf T}
$ of the odd-reduced matrices satisfying $\beta b=0$ can be interpreted as
''square roots'' of the even-reduced matrices $S$.
\end{proposition}
\begin{example}
1) Let $T^{\sqrt{S}}=\left(
\begin{array}{cc}
0 & \alpha \\
\beta & \beta \gamma
\end{array}
\right) \in {\bf T}^{\sqrt{S}}$, then $\left( T^{\sqrt{S}}\right) ^2=\alpha
\beta \left(
\begin{array}{cc}
1 & \gamma \\
0 & -1
\end{array}
\right) \in {\bf S}\setminus {\bf S}^{*}$.

2) If $\gamma =0$ in 1), then we obtain $A^2=\alpha \beta \left(
\begin{array}{cc}
1 & 0 \\
0 & -1
\end{array}
\right) \in {\bf D}\setminus {\bf D}^{*}$ (see (\ref{9}), (\ref{10}) and
compare with $D^2=\left(
\begin{array}{cc}
a^2 & 0 \\
0 & b^2
\end{array}
\right) $). This could be accepted as a definition of a square root of $%
\alpha \beta $ in some sense. Thus we have
\begin{equation}
\label{15c}
\begin{array}{ccc}
{\bf D\cdot D} & = & {\bf D}, \\ {\bf A\cdot A} & = & {\bf D},
\end{array}
\end{equation}
\noindent and the second relation could be formally considered as an ''odd
branch'' of the root $\sqrt{D}$.
\end{example}

\section{Unification of reduced supermatrices}

Now we try to unify the even- and odd-reduced matrices (\ref{4}) and (\ref{5}%
) into a common abstract object. To begin with consider the multiplication
table of all introduced sets including the even-reduced matrices products
\begin{equation}
\label{16}
\begin{array}{ccc}
{\bf S\cdot S} & = & {\bf S,} \\ {\bf D\cdot D} & = & {\bf D,} \\ {\bf %
D\cdot S} & = & {\bf S,} \\ {\bf S\cdot D} & = & {\bf S,}
\end{array}
\end{equation}

\noindent and ones for the odd-reduced matrices
\begin{equation}
\label{17}
\begin{array}{ccc}
{\bf A\cdot T} & = & {\bf S}, \\ {\bf A\cdot S} & = & {\bf T,} \\ {\bf T}%
\cdot {\bf A} & = & {\bf S}^{st}, \\ {\bf S\cdot A} & = & {\bf T}^\Pi , \\
{\bf S\cdot T} & = & {\bf S}\cup {\bf T} \\ {\bf T}\cdot {\bf S} & = & {\bf %
T.}
\end{array}
\end{equation}

\noindent  Here $st:\limfunc{{\rm Mat}} _\Lambda \left( 1|1\right)
\rightarrow \limfunc{{\rm Mat}} _\Lambda \left( 1|1\right) $ is a
supertranspose \cite{berezin}, i. e. $\left(
\begin{array}{cc}
a & \alpha \\
\beta & b
\end{array}
\right) ^{st}=\left(
\begin{array}{cc}
a & \beta \\
-\alpha & b
\end{array}
\right) $. Also we use the $\Pi $-transpose \cite{manin1} defined by $\Pi $$%
: \limfunc{{\rm Mat}} _\Lambda \left( 1|1\right) \rightarrow \limfunc{{\rm
Mat}} %
_\Lambda \left( 1|1\right) $ and
\begin{equation}
\label{pc}\left(
\begin{array}{cc}
a & \alpha \\
\beta & b
\end{array}
\right) ^\Pi =\left(
\begin{array}{cc}
b & \beta \\
\alpha & a
\end{array}
\right) .
\end{equation}
Note that the sets of matrices ${\bf S}$ and ${\bf T}$ are not closed under $%
st$ and $\Pi $ operations, but ${\bf S}^{st}\cap {\bf S=D}$ and ${\bf T}^\Pi
\cap {\bf T}={\bf A}$.

First we observe from the first two relations of (\ref{17}) that ${\bf A}$
plays a role of the left type-changing operator ${\bf A}:{\bf S}\rightarrow
{\bf T}$ and ${\bf A}:{\bf T\rightarrow S}$, while ${\bf D}$ does not change
the type. Next from the first two relations of (\ref{16}) it is obviously
seen that the sets ${\bf S}$ and ${\bf D}$ are subsemigroups. Unfortunately,
due to the next to last relation of (\ref{17}) the set ${\bf T}$ has no
clear abstract meaning. However, the last relation ${\bf T}\cdot {\bf S}=%
{\bf T}$ is important from another viewpoint: any odd-reduced morphism $%
\Lambda ^{1|1}\rightarrow \Lambda ^{1|1}$ corresponding to ${\bf T}$ can be
represented as a product of odd- and even-reduced morphisms, such that

\begin{center}
\begin{equation}
\label{d1}\setlength{\unitlength}{.15in}%
\begin{picture}(4,4) \put(0,3.1){\makebox(1,1){\large}}
\put(3.1,3.05){\makebox(1,1){\large}} \put(3.1,0){\makebox(1,1){\large}}
\put(0.3,2){\small${\bf T}$} \put(1.1,3){\vector(1,-1){1.8}}
\put(2,3.8){\small${\bf S}$} \put(1,3.5){\vector(1,0){2.3}}
\put(3.5,3.1){\vector(0,-1){1.8}}
\put(3.7,2){\small${\bf T}$} \end{picture}
\end{equation}
\end{center}

\noindent is a commutative diagram. This decomposition is crucial in the
application to the superconformal-like transformations construction (see
\cite{dup6}).

\subsection{Reduced supermatrix set semigroup}

To unify the introduced sets (\ref{16}) and (\ref{17}) we consider the
triple products
\begin{equation}
\label{18}
\begin{array}{ccc}
{\bf S}\cdot {\bf A}\cdot {\bf T} & = & {\bf S,} \\ {\bf T}\cdot {\bf A}%
\cdot {\bf T} & = & {\bf T,} \\ {\bf S}\cdot {\bf D}\cdot {\bf S} & = & {\bf %
S,} \\ {\bf T}\cdot {\bf D}\cdot {\bf S} & = & {\bf T.}
\end{array}
\end{equation}

Here we observe that the matrices ${\bf A}$ and ${\bf D}$ play the role of
''sandwich'' elements in a special ${\bf S}$ and ${\bf T}$ multiplication.
Moreover, the sandwich elements are in one-to-one correspondence with the
right sets on which they act, and so they are ''sensible from the right''.
Therefore, it is quite natural to introduce the following

\begin{definition}
{\sl A sandwich right sensible product} of the reduced supermatrix sets $%
{\bf R}={\bf S,T}$ is
\begin{equation}
\label{19}{\bf R}_1\odot {\bf R}_2\stackrel{def}{=}\left\{
\begin{array}{cc}
{\bf R}_1\cdot {\bf D\cdot R}_2, & {\bf R}_2={\bf S,} \\ {\bf R}_1\cdot {\bf %
A}\cdot {\bf R}_2, & {\bf R}_2={\bf T.}
\end{array}
\right.
\end{equation}
\end{definition}

In terms of the sandwich product instead of (\ref{18}) we obtain
\begin{equation}
\label{20}
\begin{array}{ccc}
{\bf S\odot T} & = & {\bf S,} \\ {\bf T\odot T} & = & {\bf T,} \\ {\bf %
S\odot S} & = & {\bf S,} \\ {\bf T\odot S} & = & {\bf T.}
\end{array}
\end{equation}

\begin{proposition}
The $\odot $-multiplication is associative.
\end{proposition}

\begin{proof}
Consider the relations:
\begin{equation}
\label{21}
\begin{array}{ccccc}
\left( {\bf T\odot S}\right) {\bf \odot T} & = & \left( {\bf T\cdot D\cdot S}%
\right) {\bf \cdot A\cdot T} & = & {\bf T\cdot D\cdot S\cdot A\cdot T,} \\
{\bf T\odot \left( S\odot T\right) } & = & {\bf T\cdot D\cdot \left( S\cdot
A\cdot T\right) } & = & {\bf T\cdot D\cdot S\cdot A\cdot T,}
\end{array}
\end{equation}

\noindent where the last equalities follow from the associativity of the
ordinary matrix multiplication. Therefore, $\left( {\bf T\odot S}\right)
{\bf \odot T}={\bf T\odot \left( S\odot T\right) }$. Other associativity
relations can be proved in a similar way\footnote{
We stress here that the associativity does not follow from the associativity
of the supermatrix multiplication only, but is a consequence of the special
and nontrivial set multiplication table (\ref{20}).
}.
\end{proof}
\begin{definition}
The elements ${\bf S}$ and ${\bf T}$ form a semigroup under $\odot $%
-multiplication (\ref{19}), which we call a {\sl reduced matrix set semigroup%
} and denote ${\cal RMS}_{set}$.
\end{definition}

Comparing (\ref{b1}) and (\ref{20}) we observe that the reduced matrix set
semigroup can be viewed as a right zero semigroup having two elements.
\begin{assertion}
The reduced matrix set semigroup is isomorphic to a special right zero
semigroup, i. e. ${\cal RMS}_{set}\cong {\cal Z}_R=\left\{ {\bf R=S},{\bf %
T;\odot }\right\} $.
\end{assertion}

\subsection{Scalars, anti-scalars and generalized modules}

Now we introduce the analog of $\odot $-multiplication for the reduced
matrices per se (not for sets). First we define the structure of generalized
$\Lambda $-module in $\limfunc{{\rm Hom}}_0\left( \Lambda ^{1|1},\Lambda
^{1|1}\right) $ in some alternative way, the even part of which is described
in \cite{lei1} (in the ordinary matrix theory this is a trivial fact that
the product of a matrix and a number is equal to a product of a matrix and a
diagonal matrix having this number on the diagonal).

\begin{definition}

In $\limfunc{{\rm Mat}} _\Lambda \left( 1|1\right) $ a {\sl scalar}
(matrix) $E\left( x\right) $ and {\sl anti-scalar} (matrix) ${\cal E}\left(
\chi \right) $ are defined by

\begin{equation}
\label{22}
\begin{array}{c}
E\left( x\right) \stackrel{def}{=} \left(
\begin{array}{cc}
x & 0 \\
0 & x
\end{array}
\right) \in
{\bf D}=\limfunc{{\rm Mat}} _\Lambda ^{Diag}\left( 1|1\right) ,\;x\in
\Lambda _{0,} \\ {\cal E}\left( \chi \right) \stackrel{def}{=} \left(
\begin{array}{cc}
0 & \chi \\
\chi & 0
\end{array}
\right) \in {\bf A}=\limfunc{{\rm Mat}} _\Lambda ^{Adiag}\left( 1|1\right)
,\;\chi \in \Lambda _{1.}
\end{array}
\end{equation}
\end{definition}
\begin{assertion}
The Berezin's queer subalgebra $Q_\Lambda \left( 1\right) \equiv \left(
\begin{array}{cc}
x & \chi \\
\chi & x
\end{array}
\right) \subset \limfunc{{\rm Mat}} _\Lambda \left( 1|1\right) $ \cite
{berezin} is a direct sum of the scalar and anti-scalar
\begin{equation}
\label{23}Q_\Lambda \left( 1\right) =E\left( x\right) \oplus {\cal E}\left(
\chi \right) .
\end{equation}

\end{assertion}
\begin{assertion}
The anti-scalars anticommute ${\cal E}\left( \chi _1\right) {\cal E}\left(
\chi _2\right) +{\cal E}\left( \chi _2\right) {\cal E}\left( \chi _1\right)
=0$, and so they are nilpotent.
\end{assertion}
\begin{proposition}
The structure of the generalized $\Lambda _0\oplus \Lambda _1$-module in \\
$\limfunc{{\rm Hom}}_0\left( \Lambda ^{1|1},\Lambda ^{1|1}\right) $ is
defined by action of the scalars and anti-scalars (\ref{22}).
  \end{proposition}

This means that everywhere we exchange the multiplication of supermatrices
by even and odd elements from $\Lambda $ with the multiplication by the
scalar matrices and anti-scalar ones (\ref{22}). The relations containing
the scalars are well-known \cite{lei1}, but for the anti-scalars we obtain
new dual ones. Consider their action on elements ${\bf M\in }\limfunc{{\rm
Mat}} %
_\Lambda \left( 1|1\right) $ in more detail. First we need

\begin{definition}
{\sl Left }${\cal P}$ {\sl and right }${\cal Q}$ {\sl anti-transpose }are $
\limfunc{{\rm Hom}}_0\left( \Lambda ^{1|1},\Lambda ^{1|1}\right)
\rightarrow \limfunc{{\rm Hom}}_1\left( \Lambda ^{1|1},\Lambda
^{1|1}\right) $ mappings acting on $M\in {\bf M}$ as
\begin{equation}
\label{t1}\left(
\begin{array}{cc}
a & \alpha \\
\beta & b
\end{array}
\right) ^{{\cal P}}=\left(
\begin{array}{cc}
\beta & b \\
a & \alpha
\end{array}
\right) ,
\end{equation}

\begin{equation}
\label{t2}\left(
\begin{array}{cc}
a & \alpha \\
\beta & b
\end{array}
\right) ^{{\cal Q}}=\left(
\begin{array}{cc}
\alpha & a \\
b & \beta
\end{array}
\right) .
\end{equation}
\end{definition}
\begin{corollary}
The anti-transpose is a square root of the parity changing operator (\ref{pc}%
) in the following sense
\begin{equation}
\label{atr}{\cal PQ}={\cal QP}=\Pi .
\end{equation}
\end{corollary}
\begin{assertion}
The anti-transpose satisfy
\begin{equation}
\label{atr1}
\begin{array}{ccc}
\left( {\cal E}\left( \chi \right) M\right) ^{{\cal P}} & = & \chi M \\
\left( {\cal E}\left( \chi \right) M\right) ^{{\cal Q}} & = & \chi M^\Pi \\
\left( M{\cal E}\left( \chi \right) \right) ^{{\cal P}} & = & M^\Pi \chi \\
\left( M{\cal E}\left( \chi \right) \right) ^{{\cal Q}} & = & M\chi
\end{array}
\end{equation}
\end{assertion}
Thus the concrete realization of the right, left and two-sided generalized $%
\Lambda _0\oplus \Lambda _1$-modules in $\limfunc{{\rm Hom}}_0\left(
\Lambda ^{1|1},\Lambda ^{1|1}\right) $ is determined by the actions
\begin{equation}
\label{mod1}
\begin{array}{ccc}
{\cal E}\left( \chi \right) M & = & \chi M^{
{\cal P}}, \\ M{\cal E}\left( \chi \right) & = & M^{
{\cal Q}}\chi , \\ {\cal E}\left( \chi _1\right) M{\cal E}\left( \chi
_2\right) & = & \chi _1M^\Pi \chi _2,
\end{array}
\end{equation}

\noindent together with the standard $\Lambda $-module structure \cite{lei1}
\begin{equation}
\label{mod2}
\begin{array}{ccc}
E\left( x\right) M & = & xM, \\
ME\left( x\right) & = & Mx, \\
E\left( x_1\right) ME\left( x_2\right) & = & x_1Mx_2.
\end{array}
\end{equation}

\begin{corollary}
The generalized $\Lambda _0\oplus \Lambda _1$-module relations are
\begin{equation}
\label{mod3}
\begin{array}{ccc}
\left( E\left( x\right) M\right) N & = & E\left( x\right) \left( MN\right)
\\
\left( ME\left( x\right) \right) N & = & M\left( E\left( x\right) N\right)
\\
M\left( NE\left( x\right) \right) & = & \left( MN\right) E\left( x\right) \\
\left( {\cal E}\left( \chi \right) M\right) N & = & {\cal E}\left( \chi
\right) \left( MN\right) \\ \left( M{\cal E}\left( \chi \right) \right) N &
= & M\left(
{\cal E}\left( \chi \right) N\right) \\ M\left( N{\cal E}\left( \chi \right)
\right) & = & \left( MN\right) {\cal E}\left( \chi \right)
\end{array}
\end{equation}

\noindent where $M,N\in \limfunc{{\rm Mat}} _\Lambda \left( 1|1\right) $.
\end{corollary}
\begin{proposition}
The structure of the generalized $\Lambda _0\oplus \Lambda _1$-module in \\
$\limfunc{{\rm Hom}}_1\left( \Lambda ^{1|1},\Lambda ^{1|1}\right) $ is
determined by the analogous actions of {\sl odd scalar}
\begin{equation}
\label{24o}E\left( \chi \right) \stackrel{def}{=}\left(
\begin{array}{cc}
\chi & 0 \\
0 & -\chi
\end{array}
\right) \in \limfunc{{\rm Hom}}_1\left( \Lambda ^{1|1},\Lambda
^{1|1}\right)
\end{equation}
and {\sl odd anti-scalar}
\begin{equation}
\label{25o}{\cal E}\left( x\right) \stackrel{def}{=}\left(
\begin{array}{cc}
0 & x \\
x & 0
\end{array}
\right) \in \limfunc{{\rm Hom}}_1\left( \Lambda ^{1|1},\Lambda
^{1|1}\right)
\end{equation}
respectively\footnote{
{}From (\ref{24o}) and (\ref{25o}) it is clear why we use the name ''anti-''
and not ''odd'' (as in \cite{berezin}) for the secondary diagonal matrices $%
{\cal E}\left( \chi \right) $.
}.
\end{proposition}

\subsection{Reduced supermatrix sandwich semigroup}

One way to unify the even- (\ref{4}) and odd-reduced (\ref{5}) supermatrices
into an object analogous to a semigroup is consideration of the sandwich
multiplication similar to (\ref{19}), but on the level of matrices (not
sets), by means of the scalars and anti-scalars as sandwich matrices.
Indeed, the ordinary matrix product can be written as $M_1M_2=M_1E\left(
1\right) M_2$. But we cannot find an analog of this relation using
anti-scalar, because among $\chi \in \Lambda _1$ there is no unity.
Therefore, the only possibility to include ${\cal E}\left( \chi \right) $
into equal play is consideration of sandwich elements (\ref{22}) having
arbitrary (or fixed by other special conditions) both arguments $x$ and $%
\chi $. Thus we naturally come to

\begin{definition}
{\sl A sandwich right sensible} $\Lambda _0\oplus \Lambda _1$%
{\sl -product} of the reduced supermatrices $R=S{\bf ,}T$ is
\begin{equation}
\label{26}R_1\star _XR_2\stackrel{def}{=}\left\{
\begin{array}{cc}
R_1E\left( x\right) R_2, & R_2=S
{\bf ,} \\ R_1{\cal E}\left( \chi \right) R_2, & R_2=T,
\end{array}
\right.
\end{equation}

\noindent where $X=\left\{ x,\chi \right\} \in ${\sl $\Lambda _0\oplus
\Lambda _1$}.
\end{definition}

The $\star _X$-multiplication table coincides with (\ref{20}). The
associativity can be proved similar to (\ref{21}). Therefore, we have
\begin{proposition}
Under {\sl $\Lambda _0\oplus \Lambda _1$}-multiplication the reduced
matrices form a semigroup which we call {\sl a reduced matrix sandwich
semigroup }${\cal RMSS}$.
\end{proposition}
\begin{assertion}
The reduced matrix sandwich semigroup is isomorphic to a special right zero
semigroup, i. e. ${\cal RMSS}\cong {\cal Z}_R=\left\{ R=\bigcup S\bigcup T%
{\bf ;}\star _X\right\} $.
\end{assertion}

\subsection{Direct sum of reduced supermatrices}

Another way to unify the reduced supermatrices is consideration of the
connection between them and the generalized $\Lambda _0\oplus \Lambda _1$%
-modules.

\begin{definition}
The {\sl reduced supermatrix direct space }${\cal RMDS}$ is a direct sum of
the even-reduced supermatrix space and the odd-reduced one.
\end{definition}
In terms of sets we have ${\bf R_{\oplus }=S}\oplus {\bf T}$.
\begin{assertion}
In ${\cal RMDS}$ the scalar is the Berezin's queer subalgebra $Q_\Lambda
\left( 1\right) $ (see (\ref{23})).
\end{assertion}
\begin{theorem}
In ${\cal RMDS}$ the scalars play the same role for the even-reduced
supermatrices, as the anti-scalars for the odd-reduced ones.
\end{theorem}
\begin{corollary}
The eigenvalues of even- (\ref{4}) and odd-reduced (\ref{5}) supermatrices
should be found from different equations, viz.
\begin{equation}
\label{eig1}
\begin{array}{ccc}
SV & = & E\left( x\right) V, \\
TV & = & {\cal E}\left( \chi \right) V,
\end{array}
\end{equation}

\noindent where $V$ is a column vector, and they are
\begin{equation}
\label{eig2}
\begin{array}{cc}
x_1=a, & x_2=b, \\
\chi _1=\alpha , & \chi _2=\beta .
\end{array}
\end{equation}

\noindent (see (\ref{4}) and (\ref{5})).
\end{corollary}
\begin{definition}
The characteristic functions for even- and odd-reduced supermatrices are
defined in ${\cal RMDS}$ by

\begin{equation}
\label{27}
\begin{array}{ccc}
H_S^{even}\left( x\right) & = & \limfunc{{\rm Ber}} \left( E\left( x\right)
-S\right) , \\ H_T^{odd}\left( \chi \right) & = & \limfunc{{\rm Ber}} \left(
{\cal E%
}\left( \chi \right) -T\right) .
\end{array}
\end{equation}
\end{definition}
\begin{remark}
In the standard $\Lambda $-module over $\limfunc{{\rm Mat}} _\Lambda
\left( 1|1\right) $ \cite{berezin} one derives characteristic functions and
eigenvalues for any matrix (and for odd-reduced too) from the first
equations in (\ref{eig1}) and (\ref{27}) and obtains  different results
(see, e. g. \cite{kob/nag1,urr/mor1}).
\end{remark}

Using (\ref{4}), (\ref{5}) we easily found
\begin{equation}
\label{28}
\begin{array}{ccc}
H_S^{even}\left( x\right) & = &{\displaystyle \frac{\left( x-a\right) \left(
x-b\right) }{%
\left( x-b\right) ^2}}, \\ H_T^{odd}\left( \chi \right) & = &{\displaystyle
  \frac{\left(
\chi -\alpha \right) \left( \chi -\beta \right) }{b^2}}.
\end{array}
\end{equation}

Here we observe the full symmetry between even- and odd-reduced
supermatrices (for this purpose the cancellation in the first equation was
avoided) and consistency with their $\Lambda _0\oplus \Lambda _1$%
-eigenvalues (\ref{eig2}).

The characteristic polynomial\footnote{
For a nonsupersymmetric matrix $M$ it evidently coincides with the
characteristic function $P_M\left( x\right) =H_M\left( x\right) \equiv \det
\left( Ix-M\right) $, where $I$ is a unity matrix.
}
of a supermatrix $M$ is defined by $P_M\left( M\right) =0$ and in
complicated cases is constructed from the parts of the characteristic
function $H_M\left( x\right) $ according to a special algorithm \cite
{kob/nag1,urr/mor1}. Due to existence of zero divisors in $\Lambda $ the
degree of $P_M\left( x\right) $ can be less than $n=p+q$ , $M\in \limfunc{Mat%
}_\Lambda \left( p|q\right) $. But this algorithm is not applicable for the
odd-reduced and secondary diagonal supermatrices. As before, we introduce
two dual characteristic polynomials and, using (\ref{28}), obtain the
Cayley-Hamilton theorem in ${\cal RMDS}$.
\begin{theorem}[The generalized Cayley-Hamilton theorem]
The characteristic polynomials in the reduced supermatrix direct space are
\begin{equation}
\label{29}
\begin{array}{ccc}
P_S^{even}\left( x\right) & = & \left( x-a\right) \left( x-b\right) , \\
P_T^{odd}\left( \chi \right) & = & \left( \chi -\alpha \right) \left( \chi
-\beta \right) .
\end{array}
\end{equation}

\noindent and $P_S^{even}\left( S\right) =0$ for any $S$, but $%
P_T^{odd}\left( T\right) =0$ for nilpotent $b$ only.
\end{theorem}
\begin{proof}
The even case is well-known, but for clarity we repeat it too, demonstrating
the avoiding of multiplication of a matrix by a constant and using instead
the scalars and anti-scalars (\ref{22}), i. e. the introduced $\Lambda
_0\oplus \Lambda _1$-module structure. Thus, considering simultaneously the
even and odd cases we obtain
\begin{equation}
\label{30}P_S^{even}\left( S\right) =\left( S-E\left( a\right) \right)
\left( S-E\left( b\right) \right) =\left(
\begin{array}{cc}
0 & \alpha  \\
0 & b
\end{array}
\right) \left(
\begin{array}{cc}
a & \alpha  \\
0 & 0
\end{array}
\right) =0,
\end{equation}
\begin{equation}
\label{31}
\begin{array}{c}
P_T^{odd}\left( T\right) =\left( T-
{\cal E}\left( \alpha \right) \right) \left( T-{\cal E}\left( \beta \right)
\right) = \\ \left(
\begin{array}{cc}
0 & 0 \\
\beta -\alpha  & b
\end{array}
\right) \left(
\begin{array}{cc}
0 & \alpha -\beta  \\
0 & b
\end{array}
\right) =\left(
\begin{array}{cc}
0 & 0 \\
0 & b^2
\end{array}
\right) =0.
\end{array}
\end{equation}
\end{proof}
\section{Conclusions }

We conclude that almost all above constructions are universal and ideas
mostly do not depend on size of the supermatrices under consideration. In
particular case of superconformal-like transformations it would be
interesting to use the alternative reduction introduced here in building the
objects analogous to super Riemann or semirigid surfaces, which can also lead
to
new topological-like models.

\end{document}